\documentstyle[12pt]{article}
\newcommand{\be}{\begin{equation}}
\newcommand{\ee}{\end{equation}}
\newcommand{\ba}{\begin{array}}
\newcommand{\ea}{\end{array}}
\newcommand{\ben}{\begin{enumerate}}
\newcommand{\een}{\end{enumerate}}

\newcommand{\ov}{\overline}
\newcommand{\tr}{{\rm tr}\, }
\newcommand{\Sp}{{\rm Sp}\, }

\newcommand{\Ldwa}{\,\stackrel{2}{\bigwedge}}

\newcommand{\Ph}{\mbox{$\bf Ph$}\, }

\font \msb=msbm10 scaled \magstep1

\newcommand{\rtimes}{\mbox{\msb o}\,}
\newcommand{\bR}{\mbox{\msb R} }
\newcommand{\bC}{\mbox{\msb C} }

\newcommand{\ar}{\alpha }
\newcommand{\br}{\beta }
\newcommand{\gr}{\gamma }

\newcommand{\er}{\varepsilon }

\font \eul=eufm10 scaled \magstep2

\newcommand{\gotG}{\mbox{\eul g}}
\newcommand{\gotH}{\mbox{\eul h}}

\begin{document}

\title{\bf Classical mechanical systems based on Poisson symmetry}
\author{{\bf S. Zakrzewski}  \\
\small{Department of Mathematical Methods in Physics,
University of Warsaw} \\ \small{Ho\.{z}a 74, 00-682 Warsaw, Poland} }


\date{}
\maketitle

\begin{abstract}
The existence of the theory of `twisted cotangent bundles'
(symplectic groupoids) allows to study classical mechanical
systems which are generalized in the sense that their
configurations form a Poisson manifold.  It is natural to study
from this point of view first such systems which arise in the
context of some basic physical symmetry (space-time, rotations,
etc.). We review results obtained so far in this direction.
\end{abstract}

\section{Introduction}

The theory of quantum groups \cite{D} (and Poisson groups
\cite{D:ham,S-T-S,Lu-We,Lu}) has two kind of applications in
physics:
\ben
\item {\sl hidden symmetries}. The physical model is
formulated in a standard way (mostly in the two-dimensional
`space-time') --- the formulation does not involve the
(quantum/Poisson) symmetry. This symmetry is hidden and serves
as a technical tool to establish the integrability of the model.
\item {\sl explicit symmetries}. The physical model is
drastically non-standard --- the configuration space
is considered to be
non-commutative. The quantum/Poisson symmetry controls the
non-commutativity and is also used to define the dynamics, for
instance the free one. These models are up to now only
speculative, with no direct relation to reality.
\een
In the present paper we review the results of our investigations
of classical mechanical models of the second type: with explicit
Poisson symmetries.
Our method is as follows. We pick up some typical mechanical
system whose configuration space $M$ carries a symmetry
group $G$. We `Poisson deform' $G$, i.e. choose a Poisson
structure $\pi$ on $G$ making it a Poisson group. We choose a
Poisson structure $\pi _M$ on $M$ such that the action of $G$ on
$M$ is Poisson. The symplectic groupoid of the Poisson manifold
$(M,\pi _M)$ is then considered as the {\em phase space} (it is
a natural generalization of the cotangent bundle, called also a
`twisted cotangent bundle', see
\cite{poi,poican,srni,abel} and references therein). We denote
the phase space by $\Ph (M,\pi _M)$. The action of $G$ is lifted
to $\Ph (M,\pi _M)$ (with a canonical moment map $J$ with values in
$G^*$ --- the Poisson dual of $G$) and we can consider some
dynamics compatible with the lifted action.

\section{Space-time symmetry}

Let $M$ be the Minkowski space-time and let $G$ denote the
Poincar\'{e} group acting on $M$.
$\Ph (M,\pi _M)$ is then the extended phase space of an
elementary system (without spin; the inclusion of spin requires
more subtle approach \cite{ext}) and the free dynamics is {\bf
uniquely} defined by fixing a {\em mass shell}, which is just
the inverse image of a symplectic leaf in $G^*$ by $J$. The phase
trajectories are the characteristics on the shell.

\subsection{Two-dimensional space-time}

In \cite{poi} we have considered the following Poisson structure
on the two-dimensional Minkowski space-time:
\be\label{2M}
\{x^+,x^-\}_M=\er x^+x^-,
\ee
where $x^{\pm} := x^0\pm x^1$ are the light-cone coordinates and
$\er$ --- the deformation parameter.
The (left) symplectic groupoid projection of the phase
trajectories turned out to be the constant-shape hyperboles
\be\label{hip}
 (x^+ - c^+)(x^- - c^-)= - \frac1{\er ^2 m^2}
 \ee
($m$ is the mass; the constants $c^+,c^-$ are subject to the
restriction $c^+c^-<0$). The result of the noncommutativity
(\ref{2M}) of the space-time is therefore comparable to a
presence of a repulsive force.

Taking into account that the Poisson
structure (\ref{2M}) is built of commuting vector fields:
\be\label{commut}
\pi _M = \er x^+\frac{\partial}{\partial x^+}\wedge x^-
\frac{\partial}{\partial x^-},
\ee
we obtained in \cite{poican} a realization of $\Ph (M,\pi _M)$
in the cotangent
bundle $T^*M$ (the symplectic structure remains that of the
cotangent bundle, whereas the two groupoid projections are
`deformations' of the original cotangent bundle projection, see
\cite{abel}). The cotangent bundle
projections of the phase trajectories have the following
parametric description ($p$ is the parameter):
\be
    q^+ = e^{\ar}\, \frac{\sinh \frac{\er}{2} p}
{\frac{\er}{2} m },\qquad q^- =e^{-\ar}\, \frac{\sinh
\frac{\er}{2} (p -\br )}{\frac{\er}{2} m }.
\ee
The constants $\ar,\br$ are in one-to-one correspondence
with the `scattering data' $(v_{\rm in},v_{\rm out})$:
\be
v_{\rm in}:=\lim_{p\to -\infty} \frac{q^1}{q^0} = \tan (\ar
-\frac{\er}{4}\br),\qquad
v_{\rm out}:=\lim_{p\to +\infty} \frac{q^1}{q^0} = \tan (\ar
+\frac{\er}{4}\br)
\ee
(in \cite{poican}, the factor at $\br$ was $\frac12$; we correct
that error here).
See \cite{poican,k-part} for the discussion of the role of
different projections in the symplectic groupoid and the
comparison of commuting positions $q^k$ with non-commuting
positions $x^k$.

\subsection{Four-dimensional space-time}

According to \cite{PPg}, all Poisson structures on the
Poincar\'{e} group $G$ are of coboundary type with the
$r$-matrix given by
\be\label{dec}
r= a + b + c\in
\Ldwa \gotG = \Ldwa V \oplus (V\wedge \gotH ) \oplus \Ldwa
\gotH ,
\ee
(we used the semi-direct product structure of the Poincar\'{e}
Lie algebra $\gotG=V \rtimes \gotH$ with $V$ --- translations,
$\gotH$ --- Lorentz Lie algebra), where $a$, $b$ and $c$ (the
components of $r$ coming from the decomposition on the r.h.s.
of (\ref{dec})) are listed in the table in \cite{PPg}.
For each such $r$-matrix there is exactly one Poisson structure
$\pi _M$ on $M$ such that the action is Poisson,
namely $\pi _M=r_M$, where
\be\label{rM}
r_M(x):= rx\qquad\qquad\mbox{for}\;\; x\in M
\ee
 (see \cite{poihom,standr}).
The components $a$, $b$ and $c$ correspond to the constant,
linear and quadratic part of $\pi _M$, respectively.
The construction of $\Ph (M,\pi _M)$ is relatively easy. If $r$
is triangular we can use the methods of \cite{X,standr}. If $c=0$,
the Poisson structure $\pi _M$ is linear plus constant and $\Ph
(M,\pi _M)$ is easily obtained using the cotangent bundle of an
appropriate Lie group. The only $r$ with $c\neq 0$ which is not
triangular is the case 6 of the table (cf.~correction to
\cite{PPg} in \cite{PsPg}).

\subsubsection{Triangular deformations and the abelian subcase}

If $r$ is triangular, one can realize $\Ph (M,\pi _M)$ as
(a subset of) $T^*M$ with the shifted Poisson structure \cite{standr}
\be\label{trian}
\pi _{T^*M} = r_{T^*M} + \pi _0,
\ee
where $\pi _0$ is the canonical Poisson structure of $T^*M$
(the left groupoid projection is here identified with the
original cotangent bundle projection; note that this projection
maps $\pi _{T^*M}$ on $\pi _M$).
The first term after the equality sign is defined using the
action of $G$ on $T^*M= M\times V^*$ (similarly as in
(\ref{rM})). Since the tangent action (and its dual)
involves only the action of the Lorentz part, the vertical part of
(\ref{trian}) comes only from $c$: $r_{V^*}=c_{V^*}$.
Hence the Poisson bracket between the momenta $p_k$ is purely
quadratic. The Lorentz square of the momenta, $p^2:= g^{ij}p_ip_j$,
Poisson commutes with $p_k$:
\be
\{ p^2,p_k\} = 0,
\ee
since $p^2$ is preserved by the action.
 These properties suggest
that the mass shell might be defined by $p^2 = m^2$ (cf.
\cite{Po} for a similar Ansatz on the Klein-Gordon operator).
This however should be verified by a computation of the true
mass shell (using the moment map). We leave it for the next
future.

We have made more effective calculations in the special subcase of the
triangular case: when $r$ is {\bf abelian} \cite{abel}.
In this case we can use a little different realization of $\Ph
(M,\pi _M)$ in $T^*M$ than (\ref{trian}). This time $\Ph (M,\pi
_M)$ is identified with $T^*M$ with its canonical symplectic
structure, but the left and right groupoid projections are given
by
\be\label{proj}
   M\times V^* \ni (x,p)\mapsto (x,p)_L= \exp (-\frac12 r J_0
(x,p))x\in M
\ee
\be\label{proj1}
   M\times V^* \ni (x,p)\mapsto (x,p)_R= \exp (+\frac12 r J_0
(x,p))x\in M ,
\ee
where $J_0\colon T^*M\to \gotG^*$ is the usual canonical moment
map and $r$ is treated here as a map from $\gotG^*$ to $\gotG$.
As we know from \cite{abel}, the Poisson dual of $(G,\pi )$ can be
naturally identified with $\gotG^*$ as a Poisson manifold (the
group structure of $G^*$ does not however coincide with the abelian
group structure of $\gotG^*$). Let $J\colon \Ph (M,\pi_M)\to G^*\equiv
\gotG ^*$ denote the canonical moment map
for the lift of the Poisson action of $G$ on $M$. Here is the
main result of this subsection.

\vspace{3mm}

\noindent
{\bf Theorem.} \ $J=J_0$.

\vspace{1mm}

\noindent Proof: \
We know from \cite{abel} that $\Ph (G,\pi )$ can be realized as
$T^*G$ with the left groupoid projection given by
\be
T^*G\ni \ar \mapsto \ar _L = \exp [ r(-\frac12 \ar g_{\ar }^{-1})]
g_{\ar } \exp [r(\frac12 g_{\ar }^{-1}\ar)],
\ee
where $g_{\ar}$ denotes the ordinary projection of $\ar$ on $G$.
Let $\phi $ denote the action of $G$ on $M$. We shall show that
$\Ph \phi=\Ph _0\phi$, where $\Ph \phi$ is the symplectic
groupoid lift of $\phi$ and $\Ph _0\phi$ is the ordinary
cotangent bundle lift of $\phi$ (cf.\cite{abel}). It is easy to
see that  $\Ph _0\phi $ is given by
\be
\Ph _0 \phi \; (\ar , \xi ) = g_{\ar} \xi ,
\ee
if $g^{-1}_{\ar} \ar = J_0(\xi )$ (otherwise $\Ph _0 \phi \;
(\ar , \xi ) =$ {\O}); here $\ar\in T^*G$, $\xi\in T^*M$.
The graph of $\Ph _0 \phi$ is a lagrangian submanifold of the symplectic
manifold
 $T^*M\times\overline{T^*G\times T^*M} $ (the bar
denoting the opposite symplectic structure)
which intersects the base manifold $M\times (G\times M)$ along
the graph of $\phi$. At the same time, it is contained in the
inverse image of the graph of $\phi$ by the left projection in $\Ph
(M,\pi _M)\times\overline{(\Ph (G,\pi )\times \Ph (M,\pi _M))}$,
since
$$ \phi (\ar _L , \xi _L ) =
\exp [ r(-\frac12 \ar g_{\ar }^{-1})]
g_{\ar } \exp [r(\frac12 g_{\ar }^{-1}\ar)] \cdot
\exp (-\frac12 r J_0
(\xi ))x_{\xi} =
\exp [ r(-\frac12 \ar g_{\ar }^{-1})]
g_{\ar } x_{\xi}$$
$$ =
\exp [ r(-\frac12 {\rm Ad} _{g_{\ar}} g_{\ar }^{-1}\ar)]
g_{\ar } x_{\xi} =
\exp [ r(-\frac12 J_0 (g_{\ar }\xi )]
g_{\ar } x_{\xi} = (\phi (\ar ,\xi ))_L
$$
(here $x_{\xi}$ is the ordinary projection of $\xi$)
 if $\Ph _0 \phi\; (\ar ,\xi )\neq \mbox{\O}$.
It follows that $\Ph _0 \phi$ is a solution of the problem of
characteristics which is the well known procedure to lift Poisson
maps \cite{CDW,qcp}, hence $\Ph \phi=\Ph _0\phi$. This
immediately implies the statement of the theorem.

\hfill $\Box$

{\em Corollary}. \ In the considered realization of $\Ph (M,\pi
_M)$, the mass shell is `non-deformed' and given by the standard
formula: $p^2 = m^2$. The phase trajectories are usual straight
line solutions. Their `abelian' projections coincide with
the usual trajectories of the free particle. Their groupoid
projections are easily calculated using (\ref{proj}). Note that
they are also straight lines. It follows from the fact that
$J_0(x,p)$ is constant on the phase trajectory, hence the
whole dependence of $(x,p)_L$ on the parameter comes from the
second occurrence of $x$, where the dependence is linear.

\subsubsection{Non-triangular deformations}

A representative example of a non-triangular $r$ is given by the
so called `$\kappa$-deformation' (cf.\cite{luki}; subcase 7 in
the table, with $\br =0$).
 The phase trajectories in this case has been calculated in
\cite{k-part}.
Using the fact that $\pi _M$ is a product of two commuting vector fields:
\be\label{commut1}
\pi _M = \er \frac{\partial}{\partial x^0}\wedge \sum_{k=1}^3 x^k
\frac{\partial}{\partial x^k},
\ee
we can, similarly as in (\ref{commut}), realize $\Ph (M,\pi _M)$
in the cotangent bundle and use the commuting
positions. The calculations show that
the intrinsic non-commuting positions are in this case `more physical'
than the commuting ones:
the dependence of velocity on the
momentum is monotonic in the first case and not monotonic in the
second.

\section{Free motion on Poisson $SU(2)$}

Each Poisson-Lie structure on $G=SU(2)$ is isomorphic to the
one (we denote it by $\pi $) defined by the following brackets
\be
\ba{rcl} \{\ar ,\gr\} & = & i\er \ar\gr \\
         \{\ar ,\ov{\gr}\} & = & i\er \ar\ov{\gr} \\
	 \{\gr ,\ov{\gr}\} & = & 0 \\
	 \{\ov{\ar} ,\ar\} & = & 2i\er |\gr |^2 \ea
\qquad\qquad
\left( \ba{rr} \ar & - \ov{\gr} \\
       \gr & \ov{\ar} \ea\right) =u\in SU(2),
\ee
where $\er $ is some (deformation) parameter.
According to the general rule \cite{S-T-S,Lu,WeLu2,qcp}, the phase
space of $(G,\pi )$ is given by the Manin group:
\be
\Ph (SU(2),\pi ) = SL(2,\bC ) = SU(2)\cdot SB(2) = G\cdot G^*,
\ee
where $G^*=SB(2)$ is the Poisson dual of $(G,\pi)$ composed of
upper-triangular matrices with positive diagonal elements:
\be
B= \left( \ba{rr} \rho & n \\ 0 & \rho ^{-1} \ea\right), \qquad
\rho >0, \; n\in \bC .
\ee
Each element  $A$ of $SL(2,\bC )$ has a unique (Iwasawa)
decomposition as a product of $u\in G$ and $B\in G^*$:
\be
A = \left( \ba{rr} a & b \\ c & d \ea\right) =
uB = \left( \ba{rr} \ar & - \ov{\gr} \\
       \gr & \ov{\ar} \ea\right)\cdot
       \left( \ba{rr} \rho & n \\ 0 & \rho ^{-1} \ea\right)
\ee
and $A\mapsto u$ is the left groupoid projection. The
symplectic structure on $\Ph (G,\pi )=SL(2,\bC )$ (see the
previous references for the definition), is given in terms of
coordinates by the following brackets
$$
\ba{ll}
\{ a,b \} = -i\er ab,  & \hspace{1cm}  \{ b,c \} = 2i\er ad , \\
\{ a,c \} = i\er ac ,  & \hspace{1cm} \{ b,d \} = i\er bd , \\
\{ a,d \} = 0 ,  & \hspace{1cm} \{ c,d \} = -i\er cd ,
\ea
$$
$$
\ba{lll}
\{ \ov{a},a \} = i\er ( |a|^2 + 2 |c|^2), & \hspace{6mm}
\{ \ov{b},b \} = i\er ( |b|^2 + 2|a|^2 + 2|d|^2), & \hspace{6mm}
\{ \ov{c},c \} = i\er |c|^2 ,\\
\{ \ov{b},a \} = 2i\er c\ov{d}, &
\hspace{6mm} \{ \ov{c},b \} = -i\er b\ov{c}, & \hspace{6mm}
\{ \ov{d},c \} = 0 \\
\{ \ov{c},a \} = 0, & \hspace{6mm}
\{ \ov{d},b \} = 2i\er a\ov{c},  & \hspace{6mm}
\{ \ov{d},d \} = i\er (|d|^2 + 2|c|^2) \\
\{ \ov{d},a \} = -i\er a\ov{d} . & \hspace{6mm}  & \hspace{6mm}
\ea
$$
We are now going to pick up a function on $\Ph (G,\pi )$ which
plays the role of the hamiltonian of a free motion. The natural
thing is to consider $G$-biinvariant functions on $SL(2,\bC )$
(they are in one-to-one correspondence with the Casimirs of the
Poisson structure on the `momentum space' $G^*$).
Each such a function is a function of the following basic example:
\be\label{free}
H = \frac12 \tr A^\dagger A = \frac12 \tr AA^\dagger
\ee
(this is also biinvariant in more general case $G=SU(N)$,
$SL(N,\bC ) = SU(N)\cdot SB(N)$). We take (\ref{free}) as the
hamiltonian of the free motion. Using the brackets above, we
obtain the following equations of motion:
\begin{eqnarray*}
\dot{a} & = & i\er ( aH - \ov{d} ) \\
\dot{b} & = & i\er ( bH  +\ov{c} ) \\
\dot{c} & = & i\er ( cH +\ov{b} ) \\
\dot{d} & = & i\er ( dH - \ov{a} )
\end{eqnarray*}
which can be written in the following compact form
\be\label{eqn}
\dot{A} = \{ H , A \} = i\er ( H\cdot A + Y\ov{A}Y),
\ee
where
$$  Y := \left( \ba{rr} 0 & -1 \\ 1 & 0 \ea\right) .$$
Since $H$ is a Casimir of $G^*$, the momentum $B$ is constant:
$\dot{B}=0$. Using the fact that for $u\in SU(2)$ we have $Y\ov{u}=
uY$, we obtain from(\ref{eqn})
\be
\dot{u}B= i\er ( HuB + uY\ov{B}Y),
\ee
or,
\be\label{Leg}
u^{-1} \dot{u} = i\er (H + Y\ov{B} Y B^{-1}) = {\em const}.
\ee
It follows that the trajectories in $\Ph (G,\pi )$
projected on $G$ coincide with the configurational trajectories
 of the usual free motion. Of course, the phase
description is different. The right hand side of
(\ref{Leg}) may be considered as a (deformed) Legendre
transformation (transforming momentum $B$ into velocity $u ^{-1}
\dot{u}$).

In order to make a closer comparison with the usual free motion
on $SU(2)$, let us note that the Poisson structure on the
`momentum space' $SB(2,\bC )$,
\be\label{mom}
\{ \zeta , w \} = -iw,\qquad \{ \ov{w},w\} =
i\frac{\sinh 2\er \zeta}{\er }, \qquad\qquad \mbox{where}
\;\; \rho \equiv e^{\er \zeta},\; n\equiv 2\er w
\ee
(calculated by the projection $A\mapsto B$) is isomorphic to the
usual linear Poisson structure
\be\label{lin}
\{ x-iy,x+iy\} = 2iz,\qquad \{ z, x+iy\} = -i(x+iy)
\ee
on ${\em su} (2)^*$ (see
\cite{GW} for a similar fact concerning $SU(N)$) by the
following transformation:
\be
\zeta = z,\qquad w = \frac{1}{\er}
  \sqrt{\frac{\sinh ^2 \er r - \sinh ^2 \er z}{r^2-z^2}}
(x+iy)\qquad (r^2:=x^2+y^2+z^2).
\ee
Recall that such an isomorphism defines a symplectomorphism of
$\Ph (G,\pi )$ with $T^*G$, transforming the symplectic groupoid
over $G^*$ into the usual symplectic group algebra ($T^*G$ as
the symplectic groupoid over $\gotG^*$).
The simplest Casimir of (\ref{mom}),
\be
R^2 := |w|^2 + \left(\frac{\sinh \er \zeta}{\er }\right) ^2=
\frac{H-1}{2\er ^2},
\ee
is related to $r^2$ (the natural Casimir of (\ref{lin})) as follows:
\be
\er R = \sinh \er r .
\ee
It follows that the relation between our hamiltonian $H$ and the
usual hamiltonian $h\equiv \frac12 r^2$ is given by
\be
H = \cosh 2\er r = \cosh (2\er \sqrt{2h}).
\ee
Let us make only two remarks at this point:

(i) if we took $H':=\frac12 (\frac{1}{2\er}{\rm arcosh}\, H)^2$
instead of $H$, we
would obtain a hamiltonian system on $SL(2,\bC)$ which is
symplectomorphic with the usual one of the free motion,

(ii) we conjecture that the quantum version of our model (which
likely exists) should have the quantized version $\hat{H}$ of
$H$ with the spectrum equal $\Sp \hat{H}= \cosh (2\er \sqrt{2\Sp
\hat{h}})$, where $\hat{h}$ is the quantum hamiltonian of the
free motion on $SU(2)$ (the laplacian). If one prefers to work
with a hamiltonian which tends to $h$ as $\er\to 0$, it is
reasonable to replace $H$ by
\be
H'' = \frac12 R^2 = \frac{H-1}{4\er ^2} = 
\frac12 \left(\frac{\sinh \er r}{\er }\right) ^2 =
\frac12 \left(\frac{\sinh \er \sqrt{2h}}{\er }\right) ^2 
\ee
(see also \cite{M} for an independent treatment of this
example). 

\section{Rotational symmetry}

We consider the rotation group $G:= SO(n)$ acting in the
fundamental representation: $M\equiv V = \bR ^n$.
It is not difficult to prove \cite{standr} that for any
classical $r$-matrix on $\gotG$, the bivector field $\pi
_M:=r_M$ (cf. (\ref{rM})) is Poisson (and the action of $G$ on $M$
is Poisson). In \cite{standr} we have constructed
also the (essential part of the structure of the) symplectic
groupoid of $(M,\pi _M)$ for the standard $r$-matrix given by
\be
r = i\er \sum_{\ar >0} \frac{X_\ar \wedge X_{-\ar}}{<X_\ar
,X_{-\ar}>},
\ee
where $X_\ar$ are root vectors relative to a Cartan
subalgebra, $\ar >0$ are `positive' roots and $<\cdot ,\cdot >$
is the Killing form.

Due to the complicated structure of the phase space in this case,
simple examples of hamiltonian systems having this symmetry,
like for example the deformed harmonic oscillator, are still under
investigation.


\begin{thebibliography}{99}

\bibitem{D} V. G. Drinfeld, {\em Quantum groups}, Proc. ICM,
Berkeley, 1986, vol.1, 789--820.

\bibitem{D:ham} V. G. Drinfeld, {\em Hamiltonian structures on
Lie groups, Lie bialgebras and the meaning of the classical
Yang-Baxter equations},  Soviet Math. Dokl. {\bf 27} (1983),
68--71.

\bibitem{S-T-S} M. A. Semenov-Tian-Shansky, {\em Dressing
transformations and Poisson Lie group actions},
Publ. Res. Inst. Math. Sci.,
Kyoto University {\bf 21} (1985), 1237--1260.

\bibitem{Lu-We} J.-H. Lu and A. Weinstein, {\em Poisson Lie Groups,
Dressing Transformations and Bruhat Decompositions},
J. Diff. Geom. {\bf 31} (1990), 501--526.

\bibitem{Lu} J.-H. Lu, {\em Multiplicative and affine Poisson
structures on Lie groups}, Ph.D. Thesis, University of
California, Berkeley (1990).


\bibitem{poi} S. Zakrzewski,
{\em Poisson space-time symmetry and corresponding
elementary systems},
in: ``Quantum Symmetries'', Proceedings of the II International
Wigner Symposium, Goslar 1991, H.D.~Doebner and V.K.~Dobrev (Eds.),
pp. 111--123.

\bibitem{poican} S. Zakrzewski, {\em Poisson Poincar\'{e}
particle and canonical variables}, in:
``Generalized Symmetries'',
Proceedings of the International Symposium on Mathematical
Physics, Clausthal, July 27--29, 1993, H.-D. Doebner, V.K.
Dobrev and A.G. Ushveridze (Eds.), 1994, pp. 165--171.


\bibitem{srni} S. Zakrzewski,
{\em Symplectic models of groups with noncommutative spaces}, in:
 Proceedings of the 12th Winter School on Geometry and
Topology, Srni, 11--18 January, 1992, Supplemento ai
rendiconti del Circolo Matematico di Palermo, serie II, No. 32
(1993), 185--194.


\bibitem{abel} S. Zakrzewski, {\em Geometric quantization of
Poisson groups --- diagonal and soft deformations},
Proceedings of the Taniguchi Symposium {\sl
Symplectic geometry and quantization problems}, Sanda (1993),
Y.~Maeda, H.~Omori and A.~Weinstein (Eds.), Contemporary
Mathematics {\bf 179}, 1994, 271--285.

\bibitem{ext} S.~Zakrzewski, {\em Extended phase space for a spinning
particle}, J. Phys.~A: Math. Gen. {\bf 28} (1995) 7347--7357.

\bibitem{k-part} S.~Zakrzewski, {\em On the classical $\kappa
$-particle},
in: ``Quantum Groups, Formalism and Applications'',
Proceedings of the XXX Winter School on Theoretical Physics
14--26 February 1994, Karpacz, J.~Lukierski, Z.~Popowicz,
J.~Sobczyk (eds.), Polish Scientific Publishers PWN, Warsaw
1995, pp.~573--577. Also: hep-th{/}9412098.

\bibitem{PPg} S.~Zakrzewski,
{\em Poisson Poincar\'{e} groups},
in: ``Quantum Groups, Formalism and Applications'',
Proceedings of the XXX Winter School on Theoretical Physics
14--26 February 1994, Karpacz, J.~Lukierski, Z.~Popowicz,
J.~Sobczyk (eds.), Polish Scientific Publishers PWN, Warsaw
1995, pp.~433--439. Also: hep-th{/}9412099.

\bibitem{poihom} S.~Zakrzewski, {\em Poisson homogeneous spaces},
in: ``Quantum Groups, Formalism and Applications'',
Proceedings of the XXX Winter School on Theoretical Physics
14--26 February 1994, Karpacz, J.~Lukierski, Z.~Popowicz,
J.~Sobczyk (eds.), Polish Scientific Publishers PWN, Warsaw
1995, pp.~629--639. Also: hep-th{/}9412101.

\bibitem{standr} S. Zakrzewski, {\em Phase spaces related to standard
classical $r$-matrices}, J. Phys.~A: Math. Gen. {\bf 29}
(1996) 1841--1857.  

\bibitem{X} P.  Xu, {\em Poisson manifolds associated with group
actions and classical triangular $r$-matrices},
J.~Funct.~Anal.~{\bf 112}, No.~1 (1993), 218--240.

\bibitem{PsPg} S.~Zakrzewski,
{\em Poisson structures on Poincar\'{e} group},  
Comm. Math. Phys., to appear. 


\bibitem{Po} P. Podle\'{s}, {\em Solutions of Klein-Gordon 
and Dirac Equations on Quantum Minkowski Spaces}, 
 Comm. Math. Phys. {\bf 181} (1996), 569--585.

\bibitem{CDW} A. Coste, P. Dazord and A. Weinstein,
{\em Groupo{\"\i}des symplectiques},
Publications du D\'{e}partment de Math\'{e}matiques,
Universit\'{e} Claude Bernard Lyon I (1987).

\bibitem{qcp}  S. Zakrzewski, {\em Quantum and classical
pseudogroups. Part I and II},
 Comm. Math. Phys. {\bf 134} (1990), 347--395.

\bibitem{luki} S. Zakrzewski,
{\em Quantum Poincar\'{e} group related to $\kappa$-Poincar\'{e}
algebra},  J. Phys. A: Math. Gen. {\bf 27} (1994), 2075--2082.

\bibitem{WeLu2} J.-H. Lu and A. Weinstein,
{\em Groupo{\"\i}des symplectiques doubles des groupes de
Lie-Poisson}, C.~R.~Acad.~Sci.~Paris S\'{e}r.~I~Math.~{\bf 309}
(1989), 951--954.

\bibitem{GW} V. Ginzburg and A. Weinstein, {\em Lie-Poisson
structure on some Poisson Lie groups}, J.~Amer.~Math.~Soc.~{\bf
5}, No.~2 (1992), 445--454.

\bibitem{M} G. Marmo, A. Simoni, A. Stern, {\em Poisson Lie
group symmetries for the isotropic rotator}, hep-th/9310145.


\end{thebibliography}
\end{document}